\documentclass[preprint,times,a4paper,12pt]{elsarticle}
\usepackage[latin1]{inputenc}
\usepackage[T1]{fontenc}
\usepackage{indentfirst}
\usepackage{graphicx, subcaption}
\usepackage[backref, colorlinks=true, bookmarksnumbered=true, bookmarksopen=true, bookmarksopenlevel=3, pdfstartview=FitH, linkcolor=blue, pdfmenubar=true, pdftoolbar=true, bookmarks=true,citecolor=blue, urlcolor=blue, filecolor=magenta,plainpages=false,pdfpagelabels,breaklinks]{hyperref}
\usepackage{mathrsfs, mathtools}
\usepackage{amsbsy, amsmath, amsfonts,amssymb, amsthm}

\journal{JMAA}
\theoremstyle{theorem}
\newtheorem{theorem}{Theorem}[section]

\newtheorem{proposition}[theorem]{Proposition}

\newtheorem{lemma}[theorem]{Lemma}
\newtheorem{remark}[theorem]{Remark}

\numberwithin{equation}{section}
\def\fin {\vskip 0pt \hfill \hbox{\vrule height 5pt width 5pt depth 0pt} \vskip 12pt}

\begin{document}

\begin{frontmatter}
\author[M]{Marcelo F. de Almeida\corref{cor1}}
\address[M]{Universidade Federal de Sergipe, Departamento de Matem\'{a}tica, 
Av. Marechal Rondon, s/n - Jardim Rosa Elze, CEP: 49100-000, Câmpus de São Cristóvão-SE, Brasil.}
\ead{nucaltiado@gmail.com}
\cortext[cor1]{Corresponding author}

\author[JP]{Juliana C. Precioso}
\address[JP]{Universidade Estadual Paulista ``Júlio de Mesquita Filho'', Departamento de Matem\'{a}tica, Rua Cristóvão Colombo, 2265 - Jardim Nazareth, CEP:15054-000, Câmpus de S\~{a}o Jos\'{e} do Rio Preto-SP, Brasil.}
\ead{precioso@ibilce.unesp.br}

\title{Existence and symmetries of solutions in Besov-Morrey spaces for a semilinear heat-wave type equation}

\begin{abstract}
This paper considers a semilinear integro-differential equation of Volterra type which interpolates semilinear heat and wave equations. Global existence of solutions is showed in spaces of Besov type based in Morrey spaces, namely Besov-Morrey spaces. Our initial data is larger than the previous works and our results provide a maximal existence class for semilinear interpolated heat-wave equation. Some symmetries, self-similarity and asymptotic behavior of solutions are also investigated in the framework of Besov-Morrey spaces.

\

\end{abstract}


\end{frontmatter}



\section{\protect\Large Introduction and Results}

This paper concerns with a semilinear time-fractional partial differential
equation (FPDE for short) which describes diverse physical phenomena and
mathematical models (see e.g. \cite{kilbas, Kilbas2, Povstenko}). More
precisely, in this paper we consider the semilinear integro-partial
differential equation in $\mathbb{R}^{n}$, which reads as
\begin{equation}\label{eq:fivp3}
\begin{cases}
u_{t}=\int_{0}^{t}r_{\alpha }(t-s)[P(D)u(s)+f(u(s))]ds,\;(x\in \mathbb{R}^{n}%
\text{ and }t>0) \\
u(0,x)=u_{0}(x),\;x\in \mathbb{R}^{n}
\end{cases}%
\end{equation}%
where $u(t)=u(t,x)=(u_{1}(t,x),\cdots ,u_{n}(t,x))$ with $n\geq 1$, $%
r_{\alpha }(t)=\nu t^{\alpha -1}/\Gamma (\alpha )$, $\Gamma (\alpha )$
denotes the gamma function, $P(D)=\Delta _{x}$ is the Laplacian operator on $%
x$-variable, $\nu $ denotes the Newtonian viscosity and $f:\mathbb{R}%
\rightarrow \mathbb{R}$ is a function satisfying
\begin{equation}
f(0)=0\text{ and }|f(a)-f(b)|\leq C|a-b|\left( |a|^{\rho -1}+|b|^{\rho
-1}\right) .
\end{equation}%
Above $\rho >1$ and $C$ is a positive constant independent of
$a,b\in \mathbb{R}$. Typical \linebreak examples of $f(u)$ are
given by $\gamma |u|^{\rho -1}u$ and $\gamma |u|^{\rho }$ for
$\gamma \in \{+,-\}$. These \linebreak nonlinearities yield a
scaling for (\ref{eq:fivp3}) which is fundamental for our approach
in Besov-Morrey spaces $\mathcal{N}_{p,\mu ,\infty }^{\sigma }$
(see (\ref{Besov_Morrey}), for the definition). These spaces are a
type of Besov spaces based on Morrey spaces and have been
introduced by H. Kozono and M. Yamazaki \cite{Kozo4} for analysis
of the Navier-Stokes equations. A number of authors have studied
PDEs (see \cite{Fer3, Fer4, Mazzucato, Sawano2, Xu}) and harmonic
analysis (see \cite{Sawano, Wang}) in this framework; for further
details, see \cite{Sickel1, Sickel2} and references therein. As
far as we know, the existence problem for (\ref{eq:fivp3}) in
Besov-Morrey spaces is new for the fractional case $\alpha \neq
1$. Formally, the problem (\ref{eq:fivp3}) is equivalent to (FPDE)
\begin{align}
& \partial _{t}^{\alpha }u=\nu P(D)u+f(u)\;\;\;\,\text{in }(0,\infty )\times
\mathbb{R}^{n}  \label{eq:fivp1} \\
& u_{t}(0)=0\;\text{and }\;u(0)=u_{0}\;\,\text{ in }\mathbb{R}^{n},
\label{eq:fivp2}
\end{align}
where $\mathbf{\partial }_{t}^{\alpha }u=D_{0|t}^{\alpha -1}u_{t}$, $u_{t}=%
\frac{\partial u}{\partial t}$ and $D_{0|t}^{\alpha -1}$ stands for the
Riemann-Liouville derivative of order $\alpha -1$, namely
\begin{equation*}
D_{0|t}^{\alpha -1}u=\frac{1}{\Gamma (2-\alpha )}\frac{\partial }{\partial t}%
\int_{0}^{t}\frac{u(s)}{(t-s)^{\alpha -1}}ds,\text{ for }t>0.
\end{equation*}%
Employing a Duhamel-type formula (see \cite[Proposition 2.1]{Hirata-Miao})
in (\ref{eq:fivp1})-(\ref{eq:fivp2}) (or (\ref{eq:fivp3})), formally we
obtain the integral equation
\begin{equation}
u(t)={L}_{\alpha }(t)u_{0}+B_{\alpha }(u)(t),  \label{int-mild}
\end{equation}
where
\begin{equation}
B_{\alpha }(u)(t)=\int_{0}^{t}L_{\alpha }(t-s)\left( \int_{0}^{s}r_{\alpha
-1}(s-\tau )f(u(\tau ))d\tau \right) ds  \label{nonterm0}
\end{equation}%
and $\{L_{\alpha }(t)\}_{t\geq 0}$ stands for the family of convolution
operators (or diffusion-wave operators) defined by
\begin{equation}
\widehat{L_{\alpha }(t)\varphi }(\xi )=E_{\alpha }(-t^{\alpha }|\xi |^{2})
\widehat{\varphi }(\xi ).  \label{semi-def}
\end{equation}
Throughout this paper a \textit{mild solution} for (\ref{eq:fivp1})-(\ref{eq:fivp2}) (or (\ref{eq:fivp3})) is a function $u(t,x)$ satisfying (\ref{int-mild})
and $u(t,x)\rightharpoonup u_{0}$ in $\mathcal{S}^{\prime }(\mathbb{R}^{n})$ as $t\rightarrow 0^{+}$. Actually, using Proposition \ref{conv_linear} below and Sobolev embedding (\ref{emb2}), we are going to show
this weak convergence in homogeneous Besov space $\dot{B}_{\infty ,\infty}^{2/(\rho -1)}$, see details in Lemmas \ref{linear-lem} and \ref{nonlin1}.
Here $E_{\alpha }(-t^{\alpha }|\xi |^{2})$ stands for Mittag-Leffler
function (see (\ref{mit1})) and $\,\widehat{\cdot }=\mathcal{F}\,$\ stands
for the Fourier transform. For\textbf{\ }$\alpha =1,$ the operator $L_{1}(t)=S(t)$ is the heat semigroup, because $E_{1}(-t|\xi
|^{2})=e^{-t\left\vert \xi \right\vert ^{2}}$. The kernel $k_{\alpha }$ of $%
L_{\alpha }(t)$ is the fundamental solution of (\ref{eq:fivp1}) with $%
f\equiv 0,$ namely%
\begin{equation}
{k}_{\alpha }(t,x)=\int_{\mathbb{R}^{n}}e^{ix\cdot \xi }E_{\alpha
}(-t^{\alpha }|\xi |^{2})d\xi ,  \label{fund1}
\end{equation}%
which, in one-dimensional case, reads as (see \cite{Fujita1})
\begin{equation*}
k_{\alpha }(t,x)=\frac{1}{\alpha }\int_{\mathbb{R}}\exp \{i\xi x-t|\xi |^{%
\frac{2}{\alpha }}e^{-i\frac{\gamma \pi }{2}sgn(\xi )}\}d\xi ,\;\;\left(
\gamma =2-\frac{2}{\alpha }\right) .
\end{equation*}

The FPDE (\ref{eq:fivp1})-(\ref{eq:fivp2}) interpolates two PDEs (see e.g.
\cite{Gripenberg}), namely semilinear wave ($\alpha =2$) and heat ($\alpha =1
$) equations, which have been widely investigated in the last years. These
PDEs present many differences in the theory of existence and asymptotic
behavior of solutions in scaling invariant spaces (critical spaces). In the
case $\alpha =1$, the FPDE (\ref{eq:fivp1}) is well documented in \textit{%
critical spaces}, see e.g. \cite{Kozo4}. Without making a complete list, we
mention $L^{p}$-spaces, $\text{weak-}L^{p}$ spaces, Besov spaces $\dot{B}%
_{p,\infty }^{s}$, Morrey spaces $\mathcal{M}_{p,\mu }$, Besov-Morrey spaces
$\mathcal{N}_{p,\mu ,\infty }^{s}$ and among others. However, there are few
papers dealing with FPDEs in those spaces when $1<\alpha<2$. In \cite{Hirata-Miao}, the
authors used Mihlin-Hörmander's theorem in order to establish $L^{p}\text{-}%
L^{r}$ estimates for Mittag-Leffler's family (\ref{semi-def}) and obtained
local well-posedness in a $L^{r}(\mathbb{R}^{n})$-framework. Using the
estimates of \cite{Hirata-Miao} and employing techniques of \cite{Cazenave,
Miao-Besov}, the authors of \cite{Miao2} showed the existence of
self-similar global solution with initial data $u_{0}\in \dot{B}_{p,\infty
}^{{n}/{p}-{2}/{(\rho -1)}}\cap \mathcal{E}_{q(r,p_{0}),r}$ (for $\mathcal{E}%
_{q(r,p_{0}),r}$, see Remark \ref{rem1}-(ii)). In \cite{Almeida-Ferreira},
the authors studied qualitative properties, like self-similarity,
antisymmetry and positivity, of global solutions for small initial data in
Morrey space $\mathcal{M}_{p,\lambda }$, $\lambda =n-\frac{2p}{\rho -1}$.
Now, let $P(D)f$ be the Riesz potential $(-\Delta _{x})^{{\beta }/{2}}f=%
\mathcal{F}^{-1}|\xi |^{\beta }\mathcal{F}f$ and $f(u(t))=h(x,t)|u(t)|^{\rho
-1}u(t)$. It is worth to mention the works \cite{Hakem,Kirane,Tatar} where
the authors, motivated by works of Fujita \cite{fujita_exp1, fujita_exp2},
established conditions for either blow up or global existence of weak
nonnegative solutions. It is not know if solutions of the Navier-Stokes
equations are smooth for all $t>0$, however Lions \cite{Lions} showed a
priori estimate
\begin{align}
\label{ineq-Lions}
\int_{0}^{T}\Vert D_{0|t}^{\gamma }u(t)\Vert _{\mathbf{L}^{2}(\mathbb{R}%
^{n})}dt\leq const.(J+J^{\frac{3}{2}}),\;\;J=\Vert u_{0}\Vert _{\mathbf{L}%
^{2}(\mathbb{R}^{n})}
\end{align}
where $0\leq \gamma <1/4$ and $u$ is a weak solution in $L^{2}((0,T);\mathbf{%
L}^{2}(\mathbb{R}^{n}))$ associated to the data $u_{0}\in \mathbf{L}^{2}(%
\mathbb{R}^{n})$, for $n\leq 4$. In \cite[Theorem 5.3]{Shinbrot}, Shinbrot
gave a step ahead showing (\ref{ineq-Lions}) for all dimensions $n$ and $%
0\leq \gamma <1/2$. This shows that solutions of the Navier-Stokes have more
smoothness in $t$ than at first appears. It seems that our initial data
class (see Theorem \ref{teo1}) is larger than the previous works and
contains strongly singular functions (see Remark \ref{rem1}-(iii)). For $\mu
=n-\frac{2}{\rho -1}$ and $\lambda =n-\frac{2p}{\rho -1}$, one has the
continuous inclusions
\begin{align}
\label{inv-spaces}
L^{q}\subset \text{weak-}L^{q}\subset \mathcal{M}_{p,\lambda }\subset
\mathcal{N}_{p,\mu ,\infty }^{\sigma }\;\;\text{ and }\;\;\dot{B}_{r,\infty
}^{k}\subset \mathcal{N}_{p,\mu ,\infty }^{\sigma }
\end{align}%
provided that $\frac{n}{q}=\frac{n-\lambda }{p}=-\sigma +\frac{n-\mu }{p}=-k+%
\frac{n}{r}$, where $\sigma =\frac{n-\mu }{p}-\frac{2}{\rho -1}$, $k=\frac{n%
}{r}-\frac{2}{\rho -1}$ and $1\leq q\leq r\leq p<\frac{n(\rho -1)}{2}$ (all
spaces in (\ref{inv-spaces}) are invariant by the scaling (\ref{scal2})). Moreover, there is no known existence of solutions for (\ref{eq:fivp1})-(\ref{eq:fivp2}) in a class such that  $X\varsupsetneq \mathcal{N}_{p,\mu ,\infty }^{\sigma }$. In this sense, we provide a maximal existence class for (\ref{eq:fivp1})-(\ref{eq:fivp2}) also we improve the well-posedeness result in \cite{Almeida-Ferreira,Miao2}
\bigskip

One of the aims of this paper is to establish the existence of
solutions for (\ref{eq:fivp1})-(\ref{eq:fivp2}) in the framework
of Besov-Morrey spaces. For that matter, we obtain estimates in
Sobolev-Morrey and Besov-Morrey spaces for the diffusion-wave
\linebreak operator $L_{\alpha }(t)$ (see Lemma \ref{semigroup})
which could have an interest of its own. Furthermore, some
symmetries properties, self-similarity and asymptotic behavior of
solutions are also investigated. We perform a scaling analysis in
order to choose the correct indexes of spaces such that their
norms are invariant under the scaling (\ref{scal1}). Indeed, it is
well known that if $u$ solves (\ref{eq:fivp1}) with $f(u)=\gamma
|u|^{\rho -1}u$ then, for each $\lambda >0$,
the rescaled function $u_{\lambda }(t,x)=\lambda ^{\frac{2}{\rho -1}%
}u(\lambda ^{\frac{2}{\alpha }}t,\lambda x)$ is also a solution. This leads
us to define a scaling map for (\ref{eq:fivp1}) as
\begin{equation}
u(t,x)\mapsto u_{\lambda }(t,x).  \label{scal1}
\end{equation}%
Making $t\rightarrow 0^{+}$ in (\ref{scal1}), this map induces the following
scaling for initial data
\begin{equation}
u_{0}(x)\mapsto u_{0{\lambda }}(x)=\lambda ^{\frac{2}{\rho -1}}u_{0}(\lambda
x).  \label{scal2}
\end{equation}%
Solutions invariant by (\ref{scal1}), namely $u(t,x)=u_{\lambda }(t,x)$, are
called \textit{forward self-similar solutions}.

Let $BC((0,\infty ),X)$ be the class of bounded functions from $(0,\infty )$
to a Banach space $X$. For $1<p\leq q<\infty ,$ we define our ambient space
based on Besov-Morrey type spaces (see (\ref{Besov_Morrey})) as
\begin{equation}
X_{q}^{p}=\{u\in BC((0,\infty );\mathcal{N}_{p,\mu ,\infty }^{\sigma
})\,:\,t^{\eta }u\in BC((0,\infty );\mathcal{M}_{q,\mu })\},  \label{space1}
\end{equation}%
which is a Banach space endowed with the norm
\begin{align}\label{norm1}
\Vert u\Vert _{X_{q}^{p}}:=\sup_{t>0}\Vert u(t,\cdot )\Vert _{\mathcal{N}%
_{p,\mu ,\infty }^{\sigma }}+\sup_{t>0}t^{\eta }\Vert u(t,\cdot )\Vert _{%
\mathcal{M}_{q,\mu }}.
\end{align}%
Here $\eta \in \mathbb{R}$ and $\sigma <0$ are given by
\begin{equation}
\eta =\frac{\alpha }{2}\left( \frac{2}{\rho -1}-\frac{n-\mu }{q}\right)
\text{ and }\,\sigma =\frac{n-\mu }{p}-\frac{2}{\rho -1},  \label{param1}
\end{equation}%
where these values have been chosen in such a way that the norm (\ref{norm1}) is invariant under the scaling map (\ref{scal1}).

\subsection{Main results}

In what follows, we state our main results.\smallskip

\begin{theorem}[Well-posedness]
\label{teo1} Let $n\geq 1$, $1\leq \alpha <2$, $1<\{\rho ,\,p\}\leq q<\infty
$, and $0\leq \mu <n$ be such that
 \begin{align}
 \frac{2}{\rho-1}-\frac{2}{\alpha\rho}<\frac{n-\mu}{q}<\frac{2}{\alpha(\rho-1)} \text{ and } \frac{n-\mu}{p}<\frac{2}{\rho-1}\label{H1}.
  \end{align}
\begin{itemize}
\item[\textbf{(i)}] (Existence and uniqueness) There are $\varepsilon >0$
and $\delta =\delta (\varepsilon )$ such that if $\Vert u_{0}\Vert _{%
\mathcal{N}_{p,\mu ,\infty }^{\sigma }}\leq \delta $ then the problem (\ref{eq:fivp3}) has a mild solution $u\in X_{q}^{p}$ which is unique in the
closed ball $\left\{u\in X_{q}^{p};\left\Vert u\right\Vert _{X_{q}^{p}}\leq
2\varepsilon\right\}$. Also, $u(t)\rightharpoonup u_{0}$ in the weak$%
-\ast $ topology of $\dot{B}_{\infty ,\infty }^{2/(\rho -1)}$ as $%
t\rightarrow 0^{+}$.

\item[\textbf{(ii)}] (Continuous dependence on data) Consider the ball $$
\mathcal{D}_{\delta }=\{u_{0}\in \mathcal{N}_{p,\mu ,\infty }^{\sigma
};\Vert u_{0}\Vert _{\mathcal{N}_{p,\mu ,\infty }^{\sigma }}\leq \delta \}$$
in the space $\mathcal{N}_{p,\mu ,\infty }^{\sigma }$. The data-solution map $u_{0}\in
\mathcal{D}_{\delta }\longmapsto u\in X_{q}^{p}$ is Lipschitz
continuous.
\end{itemize}
\end{theorem}

\bigskip

\begin{remark}
\label{rem1} \

\begin{itemize}
\item[\textbf{(i)}] Let $l>0$ be such that $\{p,\rho \}\leq q\leq l$ and $%
(n-\mu )/q=n/l$. By (\ref{H1}) it follows that $\alpha n(\rho -1)<2l<\alpha
n(\rho -1)\rho $ for $1\leq \alpha <2$. For every $a\in \mathcal{N}_{p,\mu
,\infty }^{\sigma }$ satisfying the assumptions of Theorem \ref{teo1}, there
exists a unique solution $u(t,x)$ of (\ref{eq:fivp1}) in $L^{\infty
}((0,\infty );\mathcal{N}_{p,\mu ,\infty }^{\sigma })$ such that $\Vert
u(t,\cdot )\Vert _{\mathcal{M}_{q,\mu }}\leq C\,t^{-\alpha /(\rho -1)+\alpha
(n-\mu )/2q}$. In particular, we recover Theorem 1 of \cite{Kozo4}.

\item[\textbf{(ii)}] Under the assumptions of Theorem \ref{teo1}, for $\mu =0
$ and $q\leq r\leq p$, we reobtain the result in \cite{Miao2}. Indeed, in
view of $\mathcal{N}_{r,0,\infty }^{k}=\dot{B}_{r,\infty }^{k}$ and
proceeding as in Lemma \ref{linear-lem} with $(p,q)=(r,q),$ one has
\begin{align*}
\Vert u_{0}\Vert _{\mathcal{E}_{q(r,p_{0}),r}}& :=\sup_{t>0}t^{\frac{1}{%
q(p_{0},r)}}\Vert L_{\alpha }(t)u_{0}\Vert _{r} \\
& =\sup_{t>0}t^{\frac{\alpha }{\rho -1}-\frac{\alpha n}{2r}}\Vert L_{\alpha
}(t)u_{0}\Vert _{r}\leq C\Vert u_{0}\Vert _{\dot{B}_{r,\infty }^{\frac{n}{r}-%
\frac{2}{\rho -1}}},
\end{align*}%
where $\frac{1}{q(p_{0},r)}=\frac{n\alpha }{2}(\frac{1}{p_{0}}-\frac{1}{r})$
and $p_{0}=\frac{n(\rho -1)}{2}$. Now, using the assumption $\Vert
u_{0}\Vert _{\dot{B}_{r,\infty }^{n/r-2/(\rho -1)}}\leq \delta $ in Theorem %
\ref{teo1}(i), one obtains Theorem 1.1 of \cite{Miao2}.

\item[\textbf{(iii)}] Let $\rho>1+2/n$ and $\lambda=n-2p/(\rho-1)$ for $p>1$%
. It follows that
\begin{align*}
\mathcal{M}_{p,\lambda}\subset \mathcal{M}_{1, n-\frac{n-\lambda}{p}}\subset
\mathcal{N}^{0}_{1,\mu,\infty} \subset\mathcal{N}^{\sigma}_{p,\mu,\infty},\;
\mu=n-\frac{2}{\rho-1}
\end{align*}
in view of (\ref{emb1}) and (\ref{emb2}). Our initial data can be taken
strictly larger than those in \cite{Almeida-Ferreira}, see \cite[Example 0.10%
]{Kozo4}.

\item[\textbf{(iv)}] Suppose that $n=1$ and let $P(D)=D_{\theta }^{\beta }$
be the Riesz-Feller operator which is given by $\widehat{D_{\theta }^{\beta
}\varphi }(\xi )=\psi _{\beta }^{\theta }(\xi )\widehat{\varphi }(\xi )$,
where $\psi _{\beta }^{\theta }(\xi )=-|\xi |^{\beta }e^{i(sgn\,\xi )\frac{%
\pi \theta }{2}}$ with $0<\beta \leq 2$ and $|\theta |\leq \min \{\beta
,2-\beta \}$, $\xi \in \mathbb{R}$. Hence (see e.g. \cite{Mainard}), the
diffusion-wave operator $L_{\alpha }(t)$ reads as
\begin{equation*}
\widehat{L_{\beta ,\alpha }^{\theta }(t)\varphi }(\xi )=E_{\alpha
}[-t^{\alpha }|\xi |^{\beta }e^{i(sgn\,\xi )\frac{\pi \theta }{2}}]\widehat{%
\varphi }(\xi )
\end{equation*}%
which has kernel
\begin{equation*}
k_{\beta ,\alpha }^{\theta }(t,x)=%
\begin{cases}
\int_{\mathbb{R}}\exp \{i\xi x-t|\xi |^{\beta }e^{-i\frac{\theta \pi }{2}%
sgn(\xi )}\}d\xi ,\;\;\alpha =1 \\
\\
\int_{\mathbb{R}}e^{i\xi x}E_{\alpha }[-t^{\alpha }\psi _{\beta }^{\theta
}(\xi )]d\xi ,\;\;\;1<\alpha <2.%
\end{cases}%
\end{equation*}%
If $(1\leq \alpha <2)$ \text{ and } $(\beta =2)$, Theorem \ref{teo1}
give us an insight on how to proceeds on the study of SFPDEs (stochastic
fractional partial differential equations)
\begin{equation*}
\frac{\partial ^{\alpha }u}{\partial t}(t,x)=D_{\theta }^{\beta
}u(t,x)+g(t,x,u(t,x))+\sum_{k=1}^{n}\frac{\partial ^{k}h_{k}}{\partial x^{k}}%
+f(t,x,u(t,x))\frac{\partial ^{2}W(t,x)}{\partial t\partial x}
\end{equation*}%
with datum $u_{0}$ in spaces more singular than $L^{p}(\mathbb{R})$ spaces.
Here, the functions $f,g,h_{k}$ satisfy Lipschitz and certain growth
conditions (see e.g. \cite{Bonaccorsi} and \cite{Niu}).
\end{itemize}
\end{remark}

\bigskip

Let $O(n)$ be the orthogonal matrix group in $\mathbb{R}^{n}$ and let $%
\mathcal{G}$ be a subset of $O(n).$ If $h(x)=h(Mx)$ and $h(x)=-h(Mx)$, for
every $M\in \mathcal{G}$, \smallskip\ then $h$ is said \textit{even} (or
symmetric) and \textit{odd} (or antisymmetric) under the action of $\mathcal{%
G}$, respectively.

\begin{theorem}
\label{teo2} Assume the hypotheses of Theorem \ref{teo1}. Let
$f(u)=\gamma |u|^{\rho -1}u.$

\begin{description}
\item[(i)] (Symmetry and antisymmetry) The solution $u(x,t)$ is
antisymmetric (resp. symmetric) for $t>0$, when $u_{0}$ is antisymmetric
(resp. symmetric) under the action of $\mathcal{G}.$

\item[(ii)] (Self-similarity) Let $u_{0}$ be a homogeneous function of
degree $-\frac{2}{\rho -1}$, then the mild solution given in Theorem \ref{teo1} is self-similar.
\end{description}
\end{theorem}

\bigskip

\begin{remark}
\label{rem2} If $\mathcal{G}=O(n)$ we have radial symmetry. Indeed, it
follows from Theorem \ref{teo2}(i) that if $u_{0}$ is radially symmetric
then $u(x,t)$ is \textbf{radially symmetric} for all $t>0$.
\end{remark}

Also, we prove an asymptotic behavior result of the solutions
obtained in \linebreak Theorem \ref{teo1} as $t\rightarrow \infty
$.

\begin{theorem}
\label{teo3}Assume the hypotheses of Theorem \ref{teo1}. Let $u$ and $v$ be
two global mild solutions for (\ref{eq:fivp3}) given by Theorem \ref{teo1},
with respective data $u_{0}$ and $v_{0}$. We have that
\begin{equation}\label{A1}
\lim_{t\rightarrow +\infty }\Vert u(\cdot ,t)-v(\cdot ,t)\Vert _{\mathcal{N}%
_{p,\mu ,\infty }^{\sigma }}=\lim_{t\rightarrow +\infty }t^{\eta }\Vert
u(\cdot ,t)-v(\cdot ,t)\Vert _{\mathcal{M}_{q,\mu }}=0
\end{equation}%
if and only if
\begin{equation}\label{A2}
\lim_{t\rightarrow +\infty }\Vert L_{\alpha }(t)(u_{0}-v_{0})\Vert _{%
\mathcal{N}_{p,\mu ,\infty }^{\sigma }}+t^{\eta }\Vert L_{\alpha
}(t)(u_{0}-v_{0})\Vert _{\mathcal{M}_{q,\mu }}=0.
\end{equation}
\end{theorem}

The manuscript is organized as follows. In Section \ref{preliminaries},
basic properties of Sobolev-Morrey, Besov-Morrey spaces and Mittag-Leffler
functions are reviewed. Section \ref{key_estimates} is devoted to estimates for operators coming from (\ref{int-mild}). Proofs of the theorems are performed in Section \ref{P_thm}.

\section{Preliminaries}

\label{preliminaries} In this section we collect some well-known properties
about Sobolev-Morrey and Besov-Morey spaces. Also, we recall properties of
the Mittag-Leffler functions.

\subsection{\protect Besov-Morrey space}

The basic properties of Morrey and Besov-Morrey spaces is reviewed in the
present subsection for the reader convenience, more details can be found in
\cite{Kato-Morrey, Kozo4, Kozo5, Peetre, Taylor}.

Let $Q_{r}(x_{0})$ be the open ball in $\mathbb{R}^{n}$ centered at $x_{0}$
and with radius $r>0$. Given two parameters $1\leq p<\infty $ and $0\leq \mu
<n$, the Morrey spaces $\mathcal{M}_{p,\mu }=\mathcal{M}_{p,\mu }(\mathbb{R}%
^{n})$ is defined to be the set of functions $f\in L^{p}(Q_{r}(x_{0}))$ such
that
\begin{equation}
\Vert f\Vert _{p,\mu }:=\sup_{x_{0}\in \mathbb{R}^{n}}\sup_{r>0}r^{-\frac{%
\mu }{p}}\Vert f\Vert _{L^{p}(Q_{r}(x_{0}))}<\infty   \label{norm-Morrey}
\end{equation}%
which is a Banach space endowed with norm (\ref{norm-Morrey}). For $s\in
\mathbb{R}$ and $1\leq p<\infty ,$ the homogeneous Sobolev-Morrey space $%
\mathcal{M}_{p,\mu }^{s}=(-\Delta )^{-s/2}\mathcal{M}_{p,\mu }$ is the
Banach space with norm
\begin{equation}
\left\Vert f\right\Vert _{\mathcal{M}_{p,\mu }^{s}}=\left\Vert (-\Delta
)^{s/2}f\right\Vert _{p,\mu }.  \label{norm-SM}
\end{equation}%
Taking $p=1$, we have $\Vert f\Vert _{L^{1}(Q_{r}(x_{0}))}$ denotes the
total variation of $f$ on open ball $Q_{r}(x_{0})$ and $\mathcal{M}_{1,\mu }$
stands for space of signed measures. In particular, $\mathcal{M}_{1,0}=%
\mathcal{M}$ is the space of finite measures. For $p>1,$ we have $\mathcal{M}%
_{p,0}=L^{p}$ and $\mathcal{M}_{p,0}^{s}=\dot{H}_{p}^{s}$ is the well known
Sobolev space. The space $L^{\infty }$ corresponds to $\mathcal{M}_{\infty
,\mu }$. Morrey and Sobolev-Morrey spaces presents the following scaling
\begin{equation}
\Vert f(\lambda \cdot )\Vert _{p,\mu }=\lambda ^{-\frac{n-\mu }{p}}\Vert
f\Vert _{p,\mu }  \label{scal-norm}
\end{equation}%
and
\begin{equation}
\left\Vert f(\lambda \cdot )\right\Vert _{\mathcal{M}_{p,\mu }^{s}}=\lambda
^{s-\frac{n-\mu }{p}}\left\Vert f\right\Vert _{\mathcal{M}_{p,\mu }^{s}}%
\text{,}  \label{scal-SM}
\end{equation}%
where the exponent $s-\frac{n-\mu }{p}$ is called scaling index and $s$ is
called regularity index. We have that
\begin{equation}
(-\Delta )^{l/2}\mathcal{M}_{p,\mu }^{s}=\mathcal{M}_{p,\mu }^{s-l}.
\label{deriv1}
\end{equation}%
Morrey spaces contain Lebesgue and $\text{weak-}L^{p}$, with the same
scaling indexes. Precisely, we have the continuous proper inclusions
\begin{equation}
L^{p}(\mathbb{R}^{n})\varsubsetneq \text{weak-}L^{p}(\mathbb{R}%
^{n})\varsubsetneq \mathcal{M}_{r,\mu }(\mathbb{R}^{n})  \label{includ1}
\end{equation}%
where $r<p$ and $\mu =n(1-r/p)$ (see e.g. \cite{Miyakawa1}). Let $\mathcal{S}%
(\mathbb{R}^{n})$ and $\mathcal{S}^{\prime }(\mathbb{R}^{n})$ be the
Schwartz space and the tempered distributions, respectively. Let $\varphi
\in \mathcal{S}(\mathbb{R}^{n})$ be nonnegative radial function such that
\begin{equation*}
supp(\varphi )\subset \{\xi \in \mathbb{R}^{n}\,;\,\frac{1}{2}<|\xi |<2\}\;%
\text{ and }\;\sum_{j=-\infty }^{\infty }\varphi _{j}(\xi )=1,\text{ for all
}\xi \neq 0
\end{equation*}%
where $\varphi _{j}(\xi )=\varphi (2^{-j}\xi )$. Let $\phi (x)=\mathcal{F}%
^{-1}(\varphi )(x)$ and $\phi _{j}(x)=\mathcal{F}^{-1}(\varphi
_{j})(x)=2^{jn}\phi (2^{j}x)$ where $\mathcal{F}^{-1}$ stands for inverse
Fourier transform. For $1\leq q<\infty $, $0\leq \mu <n$ and $s\in \mathbb{R}
$, the homogeneous Besov-Morrey space $\mathcal{N}_{q,\mu ,r}^{s}(\mathbb{R}%
^{n})$ ($\mathcal{N}_{q,\mu ,r}^{s}$ for short) is defined to be the set of $%
u\in \mathcal{S}^{\prime }(\mathbb{R}^{n})$, modulo polynomials $\mathcal{P}$%
, such that $\mathcal{F}^{-1}\varphi _{j}(\xi )\mathcal{F}u\in \mathcal{M}%
_{q,\mu }$ for all $j\in \mathbb{Z}$ and
\begin{align}
\label{Besov_Morrey}
\Vert u\Vert _{\mathcal{N}_{q,\mu ,r}^{s}}=%
\begin{cases}
\left( \sum_{j\in \mathbb{Z}}(2^{js}\Vert \phi _{j}\ast u\Vert _{q,\mu
})^{r}\right) ^{\frac{1}{r}}<\infty , & 1\leq r<\infty  \\
&  \\
\sup_{j\in \mathbb{Z}}2^{js}\Vert \phi _{j}\ast u\Vert _{q,\mu }<\infty , &
r=\infty .%
\end{cases}%
\end{align}%
The space $\mathcal{N}_{q,\mu ,r}^{s}$ is a Banach space and, in particular,
$\mathcal{N}_{q,0,r}^{s}=\dot{B}_{q,r}^{s}$ (case $\mu=0$) corresponds to the
homogeneous Besov space. We have the real-interpolation properties
\begin{equation*}
\mathcal{N}_{q,\mu ,r}^{s}=(\mathcal{M}_{q,\mu }^{s_{1}},\mathcal{M}_{q,\mu
}^{s_{2}})_{_{\theta ,r}}
\end{equation*}%
and
\begin{equation}
\label{interp_2}
\mathcal{N}_{q,\mu ,r}^{s}=(\mathcal{N}_{q,\mu ,r_{1}}^{s_{1}},\mathcal{N}%
_{q,\mu ,r_{2}}^{s_{2}})_{_{\theta ,r}},
\end{equation}%
for all $s_{1}\neq s_{2}$, $0<\theta <1$ and $s=(1-\theta )s_{1}+\theta s_{2}
$. Here $(X,Y)_{\theta ,r}$ stands for the real interpolation space between $%
X$ and $Y$ constructed via the $K_{\theta ,q}$-method. Recall that $(\cdot
,\cdot )_{\theta ,r}$ is an exact interpolation functor of exponent $\theta $
on the category of normed spaces.

In the next lemmas, we collect basic facts about Morrey spaces and
Besov-Morrey spaces (see e.g. \cite{Kato-Morrey, Kozo4, Taylor}).\smallskip

\begin{lemma}
\label{lem:2.1} Suppose that $s_{1},s_{2}\in\mathbb{R}$, $1\leq
p_1,p_2,p_3<\infty$ and $0\leq\mu_i<n$, $i=1,2,3$.

\begin{description}
\item[(i)] (Inclusion) If $\frac{n-\mu_1}{p_1}=\frac{n-\mu_2}{p_2}$ and $%
p_2\leq p_1$,
\begin{align}
\mathcal{M}_{p_1,\mu_1}\subset\mathcal{M}_{p_2,\mu_2}\;\text{ and }\;%
\mathcal{N}^{0}_{p_1,\mu_1,1}&\subset \mathcal{M}_{p_1,\mu_1}\subset
\mathcal{N}^0_{p_1,\mu_1,\infty} .  \label{emb1}
\end{align}

\item[(ii)] (Sobolev-type embedding) Let $j=1,2$ and $p_j,s_j$ be $p_2\leq
p_1$, $s_{1}\leq s_{2}$ such that $s_{2}-\frac{n-\mu_2}{p_2}=s_{1}-\frac{%
n-\mu_1}{p_1}$, we obtain
\begin{equation}
\mathcal{M}_{p_2,\mu}^{s_{2}}\subset\mathcal{M}_{p_1,\mu}^{s_{1}} ,
(\mu=\mu_1=\mu_2)  \label{sobolev-emb}
\end{equation}
and for every $1\leq r_2\leq r_1\leq\infty$, we have
\begin{align}
\mathcal{N}_{p_2,\mu_2,r_2}^{s_{2}}\subset\mathcal{N}_{p_1,%
\mu_1,r_1}^{s_{1}} \;\;\text{ and }\;\; \mathcal{N}_{p_2,\mu_2,r_2}^{s_2}%
\subset \dot{B}^{s_2-\frac{n-\mu_2}{p_2}}_{\infty,r_2} .  \label{emb2}
\end{align}

\item[\textbf{(iii)}] (Hölder inequality) Let $\;\frac{1}{p_3}=\frac{1}{p_2}+%
\frac{1}{p_1}$ and $\frac{\mu_3}{p_3}=\frac{\mu_2}{p_2}+\frac{\mu_1}{p_1}$.
If $f_j\in \mathcal{M}_{p_j,\mu_j}$ with $j=1,2$, then $f_1f_2\in\mathcal{M}%
_{p_3,\mu_3}$ and
\begin{equation}
\Vert f_1f_2\Vert_{p_3,\mu_3}\leq \Vert f_1\Vert_{p_1,\mu_1}\Vert
f_2\Vert_{p_2,\mu_2}.  \label{eq:holder}
\end{equation}
\end{description}
\end{lemma}

We finish this subsection recalling an estimate for certain multiplier
operators on $\mathcal{M}_{q,\mu}^{s}$ (see e.g. \cite{Kozo5}).\smallskip

\begin{lemma}
\label{pseudo}Let $m,s\in \mathbb{R}$ and $0\leq \mu <n$ and $P(\xi )\in
C^{[n/2]+1}(\mathbb{R}^{n}\backslash \{0\}).$ Assume that there is $A>0$
such that
\begin{equation}
\left\vert \frac{\partial ^{k}P}{\partial \xi ^{k}}(\xi )\right\vert \leq
A\left\vert \xi \right\vert ^{m-\left\vert k\right\vert },\text{ }
\label{cond-mult1}
\end{equation}%
for all $k\in (\mathbb{N}\cup \{0\})^{n}$ with $\left\vert k\right\vert \leq
\lbrack n/2]+1$ and for all $\xi \neq 0.$ Then, the multiplier operator $%
P(D)f=\mathcal{F}^{-1}P(\xi )\mathcal{F}f$ on $\mathcal{S}^{\prime }/%
\mathcal{P}$ satisfies the estimate
\begin{equation}\label{est-mult1}
\left\Vert P(D)f\right\Vert _{\mathcal{M}_{q,\mu }^{s-m}}\leq CA\left\Vert
f\right\Vert _{\mathcal{M}_{q,\mu }^{s}},\;\;(1<q<\infty )
\end{equation}
where $C>0$ is a constant independent of $f$, and the set $\mathcal{S}%
^{\prime }/\mathcal{P}$ \ consists in equivalence classes in $\mathcal{S}%
^{\prime }$ modulo polynomials with $n$ variables.
\end{lemma}

\subsection{Mittag-Leffler function}

In this part we collect some basic properties for Mittag-Leffler functions $%
E_{\alpha }(-t^{\alpha }|\xi |^{2})$ as well as the fundamental solution $%
k_{\alpha }$ (see (\ref{fund1})), further details can be obtained in \cite%
{Almeida-Ferreira, Fujita1, Hirata-Miao} and references therein.

Recall that Mittag-Leffler's function $E_{\alpha}(-t^{\alpha}|\xi|^{2})$ can
be defined via complex integral as
\begin{equation}
E_{\alpha}(-t^{\alpha}|\xi|^{2})=\frac{1}{2\pi i}\int_{\zeta}\frac{%
e^zz^{\alpha-1}}{z^{\alpha}+t^{\alpha}|\xi|^{2}}dz, \;\;(\alpha>0)
\label{mit1}
\end{equation}
where $\zeta$ is any Hankel's path on complex plan $\mathbb{C}$. The
integrand in (\ref{mit1}) has simple poles given by
\begin{equation}
a_{\alpha}(\xi)=|\xi|^{\frac{2}{\alpha}}e^{\frac{i\pi}{\alpha}%
},\;\;\;b_{\alpha}(\xi)=|\xi|^{\frac{2}{\alpha}}e^{-\frac{i\pi}{\alpha}},%
\text{ for }\xi\in\mathbb{R}^{n}.  \notag
\end{equation}

\begin{lemma}
\label{p-mittag} Let $1<\alpha<2$ and $k_{\alpha}$ be as in (\ref{fund1}).
We have that
\begin{equation}
L^{1}(\mathbb{R}^{n} )\ni E_{\alpha}(-|\xi|^{2})=\frac{1}{\alpha}%
(\exp(a_{\alpha}(\xi))+\exp(b_{\alpha }(\xi)))+l_{\alpha}(\xi) \;\; (n\geq1)
\label{Semi-prop1}
\end{equation}
where
\begin{equation*}
l_{\alpha}(\xi)=%
\begin{cases}
\frac{\sin(\alpha\pi)}{\pi}\int_{0}^{\infty}\frac{|\xi|^{2}s^{\alpha-1}e^{-s}%
}{s^{2\alpha}+2|\xi|^{2}s^{\alpha}\cos(\alpha\pi)+|\xi|^{4}}ds & \text{ if }%
\xi\neq0 \\
1-\frac{2}{\alpha}, & \text{ if }\xi=0.%
\end{cases}%
\end{equation*}
Moreover,
\begin{equation}
\frac{\partial^{k} k_{\alpha}}{\partial x_{i}^{k}}(t,x)=t^{-\frac{\alpha}{2}%
(k+n)}\frac{\partial^{k}}{\partial x_{i}^{k}}k_{\alpha}(1,t^{-\frac{\alpha}{2%
}}x), \;\;(t>0)  \label{mit-prop1}
\end{equation}
$k_{\alpha}(t,x)\geq0$, $P_{\alpha }(1,|x|)=\alpha k_{\alpha}(1,x)$ is a
probability measure.
\end{lemma}

\begin{lemma}
\label{fund-lemma} Let $1\leq\alpha<2$ and $0\leq\delta<2.$ There exists $%
A>0 $ such that
\begin{equation}
\left\vert \frac{\partial^{k}}{\partial\xi^{k}}\left[ \left\vert
\xi\right\vert ^{\delta}E_{\alpha}(-|\xi|^{2})\right] \right\vert \leq
A\left\vert \xi\right\vert ^{-\left\vert k\right\vert },\;\text{ } \xi\neq0,
\label{point1}
\end{equation}
for all $k\in(\mathbb{N}\cup\{0\})^{n}$ with $\left\vert k\right\vert
\leq\lbrack n/2]+1$.
\end{lemma}

We finish this section by recalling that $E_{\alpha }(-|\xi |^{2})$
coincides with $\sum_{k=0}^{\infty }\frac{(-|\xi |^{2})^{k}}{\Gamma (\alpha
k+1)}$, and then $E_{\alpha }(-|\xi |^{2})$ is finite, for all $\xi \in
\mathbb{R}^{n}$. Indeed, since Gamma function $\Gamma (x)$ can be expressed
as $\frac{1}{\Gamma (x)}=\frac{1}{2\pi i}\int_{\zeta }e^{z}z^{-x}dz$, for $%
x\in \mathbb{C}$, dominated convergence theorem yields
\begin{align}
E_{\alpha }(-|\xi |^{2})& =\frac{1}{2\pi i}\int_{\zeta }e^{z}z^{-1}\left[
\lim_{n\rightarrow \infty }\frac{1-(-z^{-\alpha }|\xi |^{2})^{n+1}}{%
1+z^{\alpha }|\xi |^{2}}\right] dz  \notag \\
& =\frac{1}{2\pi i}\int_{\zeta }e^{z}z^{-1}\left[ \lim_{n\rightarrow \infty
}\sum_{k=0}^{n}(-z^{-\alpha }|\xi |^{2})^{k}\right] dz  \notag \\
& =\lim_{n\rightarrow \infty }\sum_{k=0}^{n}(-|\xi |^{2})^{k}\frac{1}{2\pi i}%
\int_{\zeta }e^{z}z^{-(\alpha k+1)}dz=\sum_{k=0}^{\infty }\frac{(-|\xi
|^{2})^{k}}{\Gamma (\alpha k+1)},  \label{mittag-leffler-stand}
\end{align}%
for all $\alpha \geq 1$.

\bigskip\

\section{Key estimates}

\label{key_estimates}

The goal of this section is to derive estimates for Mittag-Leffler
convolution operators $\{L_{\alpha }(t)\}_{t\geq 0}$ on Sobolev-Morrey
spaces and Besov-Morrey spaces. Here and below the letter $C$ will denote
constants which can change from line to line.\smallskip

\begin{lemma}
\label{semigroup} Let $s,\beta\in\mathbb{R}$, $1\leq\alpha<2$, $1<p\leq
q<\infty$, $0\leq\mu<n,$ and $(\beta-s)+\frac{n-\mu}{p}-\frac{n-\mu}{q}<2$
where $\beta\geq s$.

\begin{itemize}
\item[(i)] There exists $C>0$ such that
\begin{equation}
\Vert L_{\alpha}(t)f\Vert_{\mathcal{M}^{\beta}_{q,\mu}}\leq Ct^{-\frac{\alpha%
}{2}(\beta-s)-\frac{\alpha}{2}(\frac{n-\mu}{p}-\frac{n-\mu}{q})}\Vert
f\Vert_{\mathcal{M}^{s}_{p,\mu}},  \label{semi1}
\end{equation}
for every $t>0$ and $f\in\mathcal{M}_{p,\mu}^{s}.$

\item[(ii)] Let $r\in [1,\infty]$, there exists $C>0$ such that
\begin{align}
\Vert L_{\alpha}(t)f\Vert_{\mathcal{N}^{\beta}_{q,\mu,r}}\leq Ct^{-\frac{%
\alpha}{2}(\beta-s)-\frac{\alpha}{2}(\frac{n-\mu} {p}-\frac{n-\mu}{q})}\Vert
f\Vert_{\mathcal{N}^{s}_{p,\mu,r}},  \label{semi2}
\end{align}
for every $f\in\mathcal{S}^{\prime}/\mathcal{P}$ and $t>0$.

\item[(iii)] Let $r\in [1,\infty]$ and $\beta>s$, there exists $C>0$ such
that
\begin{align}
\Vert L_{\alpha}(t)f\Vert_{\mathcal{N}^{\beta}_{q,\mu,1}}\leq Ct^{-\frac{%
\alpha}{2}(\beta-s)-\frac{\alpha}{2}(\frac{n-\mu} {p}-\frac{n-\mu}{q})}\Vert
f\Vert_{\mathcal{N}^{s}_{p,\mu,r}},  \label{semi3}
\end{align}
for every $f\in\mathcal{S}^{\prime}/\mathcal{P}$.
\end{itemize}
\end{lemma}

\noindent \textbf{Proof.} Let $\delta =(\beta -s)+l$, where $l=\frac{n-\mu }{%
p}-\frac{n-\mu }{q}$. Recalling that $\widehat{f_{\lambda }}(\xi )=\lambda
^{-n}\widehat{f}(\xi /\lambda )$ for $f_{\lambda }(x):=f(\lambda x)$, we let
$(-\Delta )^{\frac{\delta }{2}}L_{\alpha }(t)$ be the Fourier multiplier
defined as follows
\begin{align}
([(-\Delta )^{\frac{\delta }{2}}L_{\alpha }(t)]f)^{\wedge }(\xi )& =\widehat{%
h}_{\alpha }(\xi ,t)\widehat{f}(\xi )  \notag \\
& =t^{-\delta \frac{\alpha }{2}}\widehat{h}_{\alpha }(t^{\frac{\alpha }{2}%
}\xi ,1)\widehat{f}(\xi )  \notag \\
& =t^{{-\delta \frac{\alpha }{2}}}\left( h_{\alpha }(\cdot ,1)\ast
f_{t^{\alpha /2}}(\cdot )\right) _{t^{-\alpha /2}}^{\wedge }(\xi )  \notag \\
& :=t^{{-\delta \frac{\alpha }{2}}}\left( P(D)(f_{t^{\alpha /2}})\right)
_{t^{-\alpha /2}}^{\wedge }(\xi ),  \label{aux-op-1}
\end{align}%
where the symbol of $P(D)$ is $\widehat{h}_{\alpha }(\xi ,1)=|\xi |^{\delta
}E_{\alpha }(-|\xi |^{2})$. Noticing that $0\leq \delta <2$, it follows from
Lemma \ref{fund-lemma} that $P(\xi )$ shall to satisfy (\ref{cond-mult1})
with $m=0$. Using (\ref{scal-SM}) and (\ref{est-mult1}) we obtain
\begin{align}
\left\Vert \left( P(D)(f_{t^{\alpha /2}})\right) _{t^{-\alpha
/2}}\right\Vert _{\mathcal{M}_{p,\mu }^{s}}& =t^{-\frac{\alpha }{2}(s-\frac{%
n-\mu }{p})}\left\Vert P(D)(f_{t^{\alpha /2}})\right\Vert _{\mathcal{M}%
_{p,\mu }^{s}}  \notag \\
& \leq CA\,t^{-\frac{\alpha }{2}(s-\frac{n-\mu }{p})}\left\Vert f_{t^{\alpha
/2}}\right\Vert _{\mathcal{M}_{p,\mu }^{s}}  \notag \\
& =CA\,t^{-\frac{\alpha }{2}(s-\frac{n-\mu }{p})}t^{\frac{\alpha }{2}(s-%
\frac{n-\mu }{p})}\left\Vert f\right\Vert _{\mathcal{M}_{p,\mu }^{s}}  \notag
\\
& =CA\,\left\Vert f\right\Vert _{\mathcal{M}_{p,\mu }^{s}}.  \label{aux-op-2}
\end{align}%
Using (\ref{sobolev-emb}) with $(s_{1},s_{2})=(0,l)$ and (\ref{aux-op-1}),
we obtain
\begin{align}
\left\Vert L_{\alpha }(t)f\right\Vert _{\mathcal{M}_{q,\mu }^{\beta }}&
=\left\Vert (-\Delta )^{\frac{\beta }{2}}L_{\alpha }(t)f\right\Vert _{%
\mathcal{M}_{q,\mu }}  \notag \\
& \leq \left\Vert (-\Delta )^{\frac{\beta }{2}}L_{\alpha }(t)f\right\Vert _{%
\mathcal{M}_{p,\mu }^{l}},\text{ where }l=\frac{n-\mu }{p}-\frac{n-\mu }{q}
\notag \\
& =\left\Vert (-\Delta )^{\frac{\beta -s}{2}+\frac{l}{2}}L_{\alpha
}(t)f\right\Vert _{\mathcal{M}_{p,\mu }^{s}}=\;\;\left\Vert (-\Delta )^{%
\frac{\delta }{2}}L_{\alpha }(t)f\right\Vert _{\mathcal{M}_{p,\mu }^{s}}
\notag \\
& =t^{^{-\delta \frac{\alpha }{2}}}\left\Vert \left( P(D)(f_{t^{\alpha
/2}})\right) _{t^{-\alpha /2}}\right\Vert _{\mathcal{M}_{p,\mu }^{s}}
\label{x1} \\
& \leq CA\,t^{^{-\frac{\alpha }{2}(\beta -s)-\frac{\alpha }{2}(\frac{n-\mu }{%
p}-\frac{n-\mu }{q})}}\left\Vert f\right\Vert _{\mathcal{M}_{p,\mu }^{s}},
\label{x2}
\end{align}%
where (\ref{x2}) is obtained of (\ref{x1}) via inequality (\ref{aux-op-2}).
In order to obtain (\ref{semi2}) we recall the real interpolation $\mathcal{N%
}_{q,\mu ,r}^{\beta }=(\mathcal{M}_{q,\mu }^{\beta _{1}},\mathcal{M}_{q,\mu
}^{\beta _{2}})_{\theta ,r}$, $\mathcal{N}_{p,\mu ,r}^{s}=(\mathcal{M}%
_{p,\mu }^{s_{1}},\mathcal{M}_{p,\mu }^{s_{2}})_{\theta ,r}$ where $\beta
=(1-\theta )\beta _{1}+\theta \beta _{2}$, $\beta _{1}\neq \beta _{2}$ and $%
s=(1-\theta )s_{1}+\theta s_{2}$, $s_{1}\neq s_{2}$. Then, we have
\begin{equation}\label{semi_x2}
\Vert L_{\alpha }(t)f\Vert _{\mathcal{N}_{q,\mu ,r}^{\beta }}\leq
m_{0}^{1-\theta }m_{1}^{\theta }\Vert f\Vert _{\mathcal{N}_{p,\mu
,r}^{s}},\;\;0<\theta <1,
\end{equation}%
where $m_{i}=\Vert L_{\alpha }(t)f\Vert _{\mathcal{M}_{p,\mu
}^{s_{i}}\rightarrow \mathcal{M}_{q,\mu }^{\beta _{i}}}$. In view of (\ref{semi1}), we obtain
\begin{equation*}
m_{i}\leq CA\,t^{-\frac{\alpha }{2}(\beta _{i}-s_{i})-\frac{\alpha }{2}(%
\frac{n-\mu }{p}-\frac{n-\mu }{q})}
\end{equation*}%
and then, by inserting it into (\ref{semi_x2}), we get (\ref{semi2}). Now,
using (\ref{semi2}), it follows that
\begin{equation*}
\Vert L_{\alpha }(t)f\Vert _{\mathcal{N}_{q,\mu ,\infty }^{2\beta -s}}\leq
Ct^{-\alpha (\beta -s)-\frac{\alpha }{2}(\frac{n-\mu }{p}-\frac{n-\mu }{q}%
)}\Vert f\Vert _{\mathcal{N}_{p,\mu ,\infty }^{s}}
\end{equation*}%
and
\begin{equation*}
\Vert L_{\alpha }(t)f\Vert _{\mathcal{N}_{q,\mu ,\infty }^{s}}\leq Ct^{-%
\frac{\alpha }{2}(\frac{n-\mu }{p}-\frac{n-\mu }{q})}\Vert f\Vert _{\mathcal{%
N}_{p,\mu ,\infty }^{s}}.
\end{equation*}%
In view of (\ref{interp_2}) and $(2\beta -s)(1-1/2)+s(1/2)=\beta ,$ we have $%
\mathcal{N}_{q,\mu ,1}^{\beta }=(\mathcal{N}_{q,\mu ,\infty }^{2\beta -s},%
\mathcal{N}_{q,\mu ,\infty }^{s})_{1/2,1}$ which yields (\ref{semi3}). \fin

\begin{proposition}
\label{conv_linear} Let $\xi\in\mathbb{R}^n$. If $1\leq \alpha<2$, we have $%
\vert E_{\alpha}(-\vert\xi\vert^2)\vert\leq1$ and $E_{\alpha}(-\vert
t^{\alpha/2}\xi\vert^{2})\rightarrow 1$ as $t\rightarrow0^{+}$.
\end{proposition}

\noindent\textbf{Proof.} It is enough to make a proof for $1<\alpha<2$,
because the Lemma holds for $E_{1}(-t\vert\xi\vert^2)=e^{-t\vert\xi\vert^2}$%
. To this end, let $t=\vert\xi\vert^{\frac{2}{\alpha}}s^{\frac{1}{\alpha}}$
for $\xi\neq 0$ and using Lemma \ref{p-mittag} and $\vert
\exp(a_{\alpha}(\xi))\vert=\vert\exp(b_{\alpha}(\xi))\vert=\exp
(\vert\xi\vert^{2/\alpha}\cos(\pi/\alpha))\leq 1$ we obtain
\begin{eqnarray*}
\vert E_{\alpha}(-|\xi|^{2})\vert &\leq& \frac{2}{\alpha}+\vert
l_{\alpha}(\xi)\vert \\
&\leq&\frac{2}{\alpha}+\frac{\sin(\alpha\pi)}{\pi}\int_{0}^{\infty}\frac{%
e^{-\vert\xi\vert^{\frac{2}{\alpha}}s^{\frac{1}{\alpha}}}}{%
s^{2}+2s\cos(\alpha\pi)+1}ds \\
&\leq&\frac{2}{\alpha}+\frac{\sin(\alpha\pi)}{\pi}\int_{0}^{\infty}\frac{1}{%
s^{2}+2s\cos(\alpha\pi)+1}ds \\
&=&\frac{2}{\alpha} +(1-\frac{2}{\alpha})=1.
\end{eqnarray*}
Let $\Phi(t)=\frac{e^zz^{\alpha-1}}{z^{\alpha}+t^{\alpha}|\xi|^{2}}$. Note
that $\vert \Phi(t)\vert \in L^{1}(0,\infty)$ and $\Phi(t)\rightarrow
e^{z}/z $ as $t\rightarrow0^{+}$, using dominated converge theorem and
residue theorem we have
\begin{align}
E_{\alpha}(-t^{\alpha}|\xi|^{2})\rightarrow \frac{1}{2\pi i}\int_{\zeta}%
\frac{e^zz^{\alpha-1}}{z^{\alpha}}dz=\frac{1}{2\pi i}\left\{2\pi i
\,res\left(\frac{e^{z}}{z};\,z=0\right)\right\}=1,
\end{align}
because $res\left({e^{z}}/{z};\,z=0\right)=1$. \fin

\subsection{Linear estimates}

We start by recalling an elementary fixed point lemma whose proof can be
found in \cite{Fer2}.

\begin{lemma}
\label{Aux-Lem} Let $(X,\;\Vert\cdot\Vert)$ be a Banach space and $%
1<\rho<\infty$. Suppose that $B\,:\,X\rightarrow X$ satisfies $B(0)=0$ and
\begin{equation*}
\Vert B(x)-B(z)\Vert\leq K\Vert x-z\Vert(\Vert x\Vert^{\rho-1}+\Vert
z\Vert^{\rho-1})\text{.}
\end{equation*}
Let $R>0$ be the unique positive root of $2^{\rho}K\,R^{\rho-1}-1=0$.
Given $0<\varepsilon<R$ and $y\in X$ such that $\Vert y\Vert\leq\varepsilon$%
, there exists a solution $x\in X$ for the equation $x=y+B(x)$ which is the
unique one in the closed ball $\{z\in X;\left\Vert
z\right\Vert \leq2\varepsilon\}.$ Moreover, if $\Vert\bar{y}%
\Vert\leq\varepsilon$ and $\Vert\bar{x}\Vert\leq 2\varepsilon$ satisfies the
equation $\bar{x}=\bar{y}+B(\bar{x})$ then
\begin{equation}
\Vert x-\bar{x}\Vert\leq\frac{1}{1-2^{\rho}K\varepsilon^{\rho-1}}\Vert y-%
\bar{y}\Vert.  \label{d-s}
\end{equation}
\end{lemma}

The integral equation (\ref{int-mild}) has the form $x=y+B(x)$ in the
space $X=X_q^p$ where $y=L_{\alpha}(t)u_0$ and $B(x)=B_{\alpha}(u)$ is given by (\ref{semi-def}) and (\ref{nonterm0}), respectively. We invoke Lemma \ref{Aux-Lem} in our proofs, hence the estimates for linear and nonlinear part of (\ref{int-mild}) will be necessary.

\begin{lemma}
\label{linear-lem} Under the assumptions of the Theorem \ref{teo1}, there
exists $L>0$ such that
\begin{equation}
\Vert L_{\alpha}(t)u_0\Vert_{X_{q}^p}\leq L\Vert u_0\Vert_{\mathcal{N}%
_{p,\mu,\infty}^{\sigma}},  \label{linear}
\end{equation}
for all $u_0\in\mathcal{N}_{p,\mu,\infty}^{\sigma}$. Let $s=2/(\rho-1)$, if $%
u_0\in \dot{B}_{\infty,\infty}^{s}\,$ we obtain $L_{\alpha}(t)u_0%
\rightharpoonup u_0$ in the weak$-\ast$ topology of $\dot{B}%
_{\infty,\infty}^{s}$ as $t\rightarrow0^{+}$.
\end{lemma}

\noindent\textbf{Proof.} Notice that by (\ref{param1}) we obtain $\eta+\frac{\alpha}{2}\sigma=\frac{\alpha}{2}\left(\frac{n-\mu}{p}-\frac{n-\mu}{q}%
\right) $. Using (\ref{H1}) one has $\frac{n-\mu}{p}-\frac{n-\mu}{q}-\sigma=\frac{2}{\rho-1}-\frac{n-\mu}{q}<\frac{2}{\alpha \rho}<2$ and $\sigma<0$
which by (\ref{semi2}) and afterwards by (\ref{emb1}) and (\ref{semi3}),
respectively, give us
\begin{align}
\sup_{t>0}\Vert L_{\alpha}(t)u_0\Vert_{\mathcal{N}_{p,\mu,\infty}^{\sigma}}
+\sup_{t>0}t^{\eta}\Vert &L_{\alpha}(t)u_0\Vert_{\mathcal{M}_{q,\mu}}\leq
C\Vert a\Vert_{\mathcal{N}_{p,\mu,\infty}^{\sigma}} +\sup_{t>0}t^{\eta}\Vert
L_{\alpha}(t)u_0\Vert_{\mathcal{N}_{q,\mu,1}^{0}}  \notag \\
&\leq C\Vert u_0\Vert_{\mathcal{N}_{p,\mu,\infty}^{\sigma}} +
C\sup_{t>0}t^{\eta+\frac{\alpha}{2}\sigma-\frac{\alpha}{2}(\frac{n-\mu}{p}-%
\frac{n-\mu}{q})}\Vert u_0\Vert_{\mathcal{N}_{p,\mu,\infty}^{\sigma}}  \notag
\\
&\leq L\Vert u_0\Vert_{\mathcal{N}_{p,\mu,\infty}^{\sigma}} ,  \notag
\end{align}
this yield (\ref{linear}). It remains to check the weak$-\ast$ convergence.
To this end, let $v\in \dot{B}^{-s}_{1,1}$ the predual space of $\dot{B}%
^{s}_{\infty,\infty}$. Using Proposition \ref{conv_linear} we have
\begin{align}
\Vert L_{\alpha}(t)v - v\Vert_{\dot{B}^{-s}_{1,1}
}=\sum_{j=-\infty}^{\infty}\{2^{-js}\Vert \mathcal{F}^{-1}
\varphi_j[E_{\alpha}(-t^{\alpha}\vert\xi\vert^2)-1]\mathcal{F}v\Vert_{L^{1}(%
\mathbb{R}^n)}\} \rightarrow 0  \notag
\end{align}
as $t\rightarrow0^{+}$. Thanks to (\ref{mit-prop1}) and (\ref{mittag-leffler-stand}) one has $E_{\alpha}(-t^{\alpha}\vert\xi\vert^2)=t^{-%
\frac{\alpha}{2}n}E_{\alpha}(-\vert t^{-\frac{\alpha}{2}}\xi\vert^2)\in\mathbb{R}$ for all $t>0$ and $\xi\in\mathbb{R}^n$, it follows that
\begin{align}
\vert \langle L_{\alpha}(t)u_0 - u_0, v\rangle \vert =\vert \langle u_0,
L_{\alpha}(t)v - v\rangle \vert\leq \Vert u_0\Vert_{\dot{B}%
^{s}_{\infty,\infty}}\Vert L_{\alpha}(t)v - v\Vert_{\dot{B}%
^{-s}_{1,1}}\rightarrow 0,
\end{align}
as $t\rightarrow0^{+}$. \fin

\subsection{Nonlinear estimates}

Recall the nonlinear term in (\ref{int-mild})
\begin{equation}
B_{\alpha }(u)(t)=\int_{0}^{t}L_{\alpha }(t-s)\int_{0}^{s}r_{\alpha
-1}(s-\tau )f(u(\tau ))d\tau ds.  \label{nonterm1}
\end{equation}

\begin{lemma}[Nonlinear estimate]
\label{nonlin1} Assume the hypotheses of Theorem \ref{teo1}. There is a
constant $K>0$ such that
\begin{equation}\label{bil1}
\Vert B_{\alpha }(u)-B_{\alpha }(v)\Vert _{X_{q}^{p}}\leq K\Vert u-v\Vert
_{X_{q}^{p}}(\Vert u\Vert _{X_{q}^{p}}^{\rho -1}+\Vert v\Vert
_{X_{q}^{p}}^{\rho -1}),
\end{equation}%
for all $u,v\in X_{q}^{p}$. Moreover, we have $B_{\alpha
}(u)(t)\rightharpoonup 0$ in the weak$-\ast $ topology of $\dot{B}_{\infty
,\infty }^{2/(\rho -1)}$ as $t\rightarrow 0^{+}$.
\end{lemma}

\noindent\textbf{Proof.} The proof is divided in three steps.

\noindent\textbf{First step.} Let $\tilde{s}\in\mathbb{R}$ be such that $\sigma-\frac{n-\mu}{p}=\tilde{s}-\frac{n-\mu}{q/\rho}$. In view of (\ref{H1}) and $\alpha\geq1$ we have $\frac{n-\mu}{q/\rho}\geq \frac{2}{\rho-1}>\frac{%
n-\mu}{p}$, it follows that $\sigma<0\leq \tilde{s}$ and $p\geq q/\rho$.
Applying (\ref{sobolev-emb}) and (\ref{semi2}) afterwards (\ref{emb1}) and (\ref{eq:holder}), respectively, we have
\begin{align}
\Vert B_{\alpha}(u)(t)-&B_{\alpha}(v)(t)\Vert_{\mathcal{N}%
^{\sigma}_{p,\mu,\infty}}\leq \Vert B_{\alpha}(u)(t)-B_{\alpha}(v)(t)\Vert_{%
\mathcal{N}^{\tilde{s}}_{q/\rho,\mu,\infty}}  \notag \\
&\leq C\int_{0}^{t}(t-s)^{\gamma_1}\int_{0}^{s}r_{\alpha-1}(s-\tau)\Vert
f(u(\tau))-f(v(\tau))\Vert_{\mathcal{N}_{q/\rho,\mu,\infty}^{0}}d\tau ds
\notag \\
& \leq C\int_{0}^{t}(t-s)^{\gamma_1}\int_{0}^{s}r_{\alpha-1}(s-\tau)\Vert
f(u(\tau))-f(v(\tau))\Vert_{\mathcal{M}_{q/\rho,\mu}}d\tau ds  \notag \\
& \leq C\int_{0}^{t}(t-s)^{\gamma_1}\int_{0}^{s}r_{\alpha-1}(s-\tau)\Vert
u-v\Vert_{\mathcal{M}_{q,\mu}}(\Vert u\Vert_{\mathcal{M}_{q,\mu}}^{\rho-1}+%
\Vert v\Vert_{\mathcal{M}_{q,\mu}}^{\rho-1})d\tau ds  \label{nonlin_B1-est11}
\\
&:=\psi_{1}(t)\sup_{t>0}t^{\eta}\Vert u(t)-v(t)\Vert_{\mathcal{M}_{q,\mu}}
\sup_{t>0}t^{\eta(\rho-1)}(\Vert u(t)\Vert_{\mathcal{M}_{q,\mu}}^{\rho-1}+%
\Vert v(t)\Vert_{\mathcal{M}_{q,\mu}}^{\rho-1})  \label{nonlin_B1-est1}
\end{align}

where $r_{\alpha}(s)=s^{\alpha-1}/\Gamma(\alpha)$, $f(u(\tau))=\vert u(\tau
)\vert^{\rho-1}u(\tau)$, $\gamma_1=-\frac{\alpha}{2}\tilde{s}=-\frac{\alpha}{%
2}\sigma-\frac{\alpha}{2}(\frac{n-\mu}{q/\rho}-\frac {n-\mu}{p})$ and
\begin{align}
\psi_1(t)&=C\int_{0}^{t}(t-s)^{\gamma_1}\int_{0}^{s}r_{\alpha-1}(s-\tau)%
\tau^{-\eta\rho}d\tau ds  \notag \\
&=C\beta(1-\eta\rho,\alpha-1)\beta(\alpha-\eta,\gamma_1+1),  \label{aux_psi1}
\end{align}
for $\beta(x,y)=\int_0^1t^{x-1}(1-t)^{y-1}dt$ if $x,y>0$. Indeed, by change
of variables $\tau=zs$ and $s=t\omega$, respectively, we get
\begin{align}
\int_{0}^{t}(t-s)^{\gamma_1}\int_{0}^{s}r_{\alpha-1}(s-\tau)&\tau^{-\eta%
\rho}d\tau ds
=\int_{0}^{t}(t-s)^{\gamma_1}s^{\alpha-1-\eta\rho}\left(\int_{0}^{1}(1-z)^{%
\alpha-2}z^{-\eta\rho}dz\right)ds  \notag \\
&=\beta(1-\eta\rho,\alpha-1)t^{\alpha+\gamma_1-\eta\rho}\int_{0}^{1}(1-%
\omega)^{\gamma_1}\omega^{\alpha-1-\eta\rho}d\omega,  \notag \\
&=\beta(1-\eta\rho,\alpha-1)\beta(\alpha-\eta\rho,\gamma_1+1),
\label{betas1}
\end{align}
because by $\gamma_1=-\frac{\alpha}{2}\sigma-\frac{\alpha}{2}(\frac{n-\mu}{%
q/\rho}-\frac {n-\mu}{p})$ and (\ref{param1}) we have
\begin{align*}
\alpha+\gamma_1-\eta\rho &=\alpha+\frac{\alpha}{2}(\frac{2}{\rho-1}-\frac{%
n-\mu}{q}\rho) -\frac{\alpha}{2}(\frac{2}{\rho-1}-\frac{n-\mu}{q})\rho \\
&=\alpha +\frac{\alpha}{\rho -1} -\frac{\alpha\rho}{\rho -1}=0.
\end{align*}
Inserting (\ref{aux_psi1}) into (\ref{nonlin_B1-est1}) yields {\footnotesize
\begin{equation}
\Vert B_{\alpha}(u)(t)-B_{\alpha}(v)(t)\Vert_{\mathcal{N}^{\sigma}_{p,\mu,%
\infty}}\leq K_1\sup_{t>0}t^{\eta}\Vert u(t)-v(t)\Vert_{\mathcal{M}%
_{q,\mu}}\sup_{t>0}t^{\eta(\rho-1)}(\Vert u(t)\Vert_{\mathcal{M}%
_{q,\mu}}^{\rho-1}+\Vert v(t)\Vert_{\mathcal{M}_{q,\mu}}^{\rho-1}).
\label{aux-non1}
\end{equation}
} \bigskip \textbf{Second step.} Let $\beta=s=0$. By estimate (\ref{semi1})
and Hölder inequality (\ref{eq:holder}) we obtain
\begin{align}
\Vert B_{\alpha}(u)(t)-B_{\alpha}(v)(t)\Vert_{\mathcal{M}_{q,\mu}} & \leq
C\int_{0}^{t}(t-s)^{\gamma_{2}}\theta(s)ds  \label{nonlin-est1}
\end{align}
where $\gamma_{2}=-\frac{\alpha}{2}(\frac {n-\mu}{q/\rho}-\frac{n-\mu}{q})$
and $\theta(s)$ is given by
\begin{align}
\theta(s) & =\int_{0}^{s}r_{\alpha-1}(s-\tau)\Vert u(\tau)-v(\tau )\Vert_{%
\mathcal{M}_{q,\mu}}(\Vert u(\tau)\Vert_{\mathcal{M}_{q,\mu}}^{\rho-1}+\Vert
v(\tau)\Vert_{\mathcal{M}_{q,\mu}}^{\rho-1})d\tau .  \notag
\end{align}
Mimicking the \textit{First step} we get {\footnotesize
\begin{equation}
\Vert B_{\alpha}(u)(t)-B_{\alpha}(v)(t)\Vert_{\mathcal{M}_{q,\mu}}\leq
\psi_2(t)\sup_{t>0}t^{\eta}\Vert u(t)-v(t)\Vert_{\mathcal{M}%
_{q,\mu}}\sup_{t>0}t^{\eta(\rho-1)}(\Vert u(t)\Vert_{\mathcal{M}%
_{q,\mu}}^{\rho-1}+\Vert v(t)\Vert_{\mathcal{M}_{q,\mu}}^{\rho-1})
\label{aux-non4}
\end{equation}
} where $\psi_2(t)$ can be estimated as
\begin{align}
\psi_{2}(t)&\leq
C\beta(1-\eta\rho,\alpha-1)\beta(\alpha-\eta\rho,\gamma_2+1)t^{\alpha+%
\gamma_2-\eta\rho}=K_2 t^{-\eta},  \label{betas2}
\end{align}
because in view of (\ref{param1}) we have
\begin{equation*}
\alpha+\gamma_2-\eta\rho=\alpha -\frac{\alpha}{2}\frac{n-\mu}{q}-\frac{\alpha%
}{2}\frac{2\rho}{\rho-1}=\alpha+\frac{\alpha}{\rho-1}-\frac{\alpha\rho}{%
\rho-1}-\eta =-\eta.
\end{equation*}
Inserting (\ref{betas2}) into (\ref{aux-non4}) it follows that
{\footnotesize
\begin{equation}
t^{\eta}\Vert B_{\alpha}(u)(t)-B_{\alpha}(v)(t)\Vert_{\mathcal{M}%
_{q,\mu}}\leq K_2\sup_{t>0}t^{\eta}\Vert u(t)-v(t)\Vert_{\mathcal{M}%
_{q,\mu}}\sup_{t>0}t^{\eta(\rho-1)}(\Vert u(t)\Vert_{\mathcal{M}%
_{q,\mu}}^{\rho-1}+\Vert v(t)\Vert_{\mathcal{M}_{q,\mu}}^{\rho-1}).
\label{aux-non2}
\end{equation}
} The convergence of the beta functions appearing in (\ref{betas1}) and (\ref{betas2}) is obtained by restrictions (\ref{H1}) and $\alpha\geq 1$, because
this yields in $\gamma_1,\gamma_2>-1$ and $\eta\rho<1\leq\alpha$. It follows
that $(\frac{n-\mu}{q/\rho}-\frac{n-\mu}{q})<\frac{2}{\alpha}\leq2$ which we
have used in \textit{Second step}. Recalling (\ref{norm1}) and using (\ref{aux-non1}) and (\ref{aux-non2}) we obtain (\ref{bil1}) with $K=K_{1}+K_{2}$.

\bigskip \noindent\textbf{Third step.} As $\mathcal{S}(\mathbb{R}^{n})$ is
dense in $\dot{B}^{{2}/{(\rho-1)}}_{1,1}$ (see \cite[p. 48]{Triebel}) the
weak$-\ast$ convergence can be obtained by estimate
\begin{align}
\vert\langle B_{\alpha}(u)(t),v\rangle \vert &\leq \vert\langle
B_{\alpha}(u)(t),v-\varphi\rangle \vert + \vert\langle
B_{\alpha}(u)(t),\varphi\rangle \vert  \notag \\
&\leq \Vert B_{\alpha}(u)(t)\Vert_{\dot{B}^{\sigma-{(n-\mu)}/{p}%
}_{\infty,\infty}} \Vert v-\varphi\Vert_{\dot{B}^{{2}/{(\rho-1)}}_{1,1}} +
\vert \langle B_{\alpha}(u)(t), \varphi\rangle\vert  \notag \\
&\leq C\Vert u\Vert_{X_{q}^p}\varepsilon +C\Vert u\Vert_{X_{q}^p}^{\rho}
\Vert \varphi\Vert_{\dot{B}_{1,1}^{{2}/{(\rho-1)}}}t^{\alpha}\leq C\Vert
u\Vert_{X_{q}^p}\varepsilon \text{ as } t\rightarrow0^+,
\end{align}
because for $v\in\dot{B}^{\sigma-(n-\mu)/{p}}_{1,1}=\dot{B}^{{2}/{(\rho-1)}%
}_{1,1}$ one has $\Vert v-\varphi\Vert_{\dot{B}^{{2}/(\rho-1)}_{1,1}} \leq
\varepsilon$, for all $\varepsilon>0$. Moreover, by embedding $\mathcal{N}%
^{s}_{l,\mu,\infty}\subset\dot{B}^{s-{(n-\mu)}/{l}}_{\infty,\infty}$ (see (%
\ref{emb2})) and $\frac{n-\mu}{p}<\frac{2}{\rho-1}<\frac{n-\mu}{q/\rho}$ we have that
\begin{align}
\vert \langle B_{\alpha}(u)(t), \varphi\rangle\vert &=
\left\vert\int_{0}^t\int_{0}^sr_{\alpha-1}(s-\tau)\left\langle f(u(\tau)),
L_{\alpha}(t-s)\varphi\right\rangle d\tau ds\right\vert  \notag \\
&\leq C\int_{0}^t\int_{0}^sr_{\alpha-1}(s-\tau) \Vert\, \vert
u(\tau)\vert^{\rho-1}u(\tau)\Vert_{\dot{B}_{\infty,\infty}^{-\frac{n-\mu}{%
q/\rho}}} \Vert L_{\alpha}(t-s)\varphi\Vert_{\dot{B}_{1,1}^{\frac{n-\mu}{%
q/\rho}}} d\tau ds  \notag \\
&\leq C\int_{0}^t\int_{0}^s(s-\tau)^{\alpha-2} (t-s)^{-\frac{\alpha}{2}(%
\frac{n-\mu}{q/\rho}-\frac{2}{\rho-1})}\Vert\, \vert
u(\tau)\vert^{\rho-1}u(\tau)\Vert_{\mathcal{N}_{q/\rho,\mu,\infty}^{0}}
\Vert \varphi\Vert_{\dot{B}_{1,1}^{\frac{2}{\rho-1}}} d\tau ds  \notag \\
&\leq C\int_{0}^t\int_{0}^s(s-\tau)^{\alpha-2} (t-s)^{-\frac{\alpha}{2}(%
\frac{n-\mu}{q/\rho}-\frac{2}{\rho-1})} \tau^{-\eta\rho}d\tau ds\,\Vert
u\Vert_{X_q^p}^{\rho}\Vert \varphi\Vert_{\dot{B}_{1,1}^{\frac{2}{\rho-1}}}
\notag \\
&\leq Ct^{\alpha}\Vert u\Vert_{X_q^p}^{\rho}\Vert \varphi\Vert_{\dot{B}%
_{1,1}^{\frac{2}{\rho-1}}}\rightarrow 0\; \text{ as }\; t\rightarrow 0^+
\end{align}
and
\begin{align}
\Vert B_{\alpha}(u)(t)\Vert_{\dot{B}^{\sigma-\frac{n-\mu}{p}%
}_{\infty,\infty}} &\leq C\int_{0}^t\int_{0}^sr_{\alpha-1}(s-\tau) \Vert
L_{\alpha}(t-s)\,\vert u(\tau)\vert^{\rho-1}u(\tau)\Vert_{\dot{B}^{\sigma-%
\frac{n-\mu}{p}}_{\infty,\infty}}d\tau ds  \notag \\
&\leq \int_{0}^t\int_{0}^sr_{\alpha-1}(s-\tau) \Vert L_{\alpha}(t-s)\,\vert
u(\tau)\vert^{\rho-1}u(\tau)\Vert_{\mathcal{N}^{\sigma}_{p,\mu,\infty}}d\tau
ds  \notag \\
&\leq \int_{0}^t\int_{0}^sr_{\alpha-1}(s-\tau) \Vert L_{\alpha}(t-s)\,\vert
u(\tau)\vert^{\rho-1}u(\tau)\Vert_{\mathcal{N}^{\tilde{s}}_{q/\rho,\mu,%
\infty}}d\tau ds  \notag \\
&\leq C\int_{0}^t\int_{0}^sr_{\alpha-1}(s-\tau) (t-s)^{\gamma_1}\Vert\,\vert
u(\tau)\vert^{\rho-1}u(\tau)\Vert_{\mathcal{N}^{0}_{q/\mu,\mu,\infty}}d\tau
ds  \notag \\
&\leq C\int_{0}^t\int_{0}^sr_{\alpha-1}(s-\tau)
(t-s)^{\gamma_1}\tau^{-\eta\rho}d\tau ds\,\Vert u\Vert_{X_{q}^p}^{\rho}
\notag \\
&\leq C\,\Vert u\Vert_{X_{q}^p}^{\rho},
\end{align}
as required, this finish our proof. \fin

\section{Proof of theorems}

\label{P_thm}

\subsection{Proof of Theorem \protect\ref{teo1}}

Let $0<\varepsilon< R=\left(1/2^{\rho}K\right)^{\rho-1}$, where $K>0$ and $%
L>0$ are the constants obtained in Lemma \ref{linear-lem} and Lemma \ref{nonlin1}, respectively. Let $\delta=\varepsilon/L$, the Lemma \ref{Aux-Lem}
with $X=X_q^p$ and $y=L_{\alpha}(t)u_0$ yields the existence of an unique
global mild solution $u\in X_q$ such that $\Vert u\Vert_{X_q}\leq
\varepsilon $. Moreover, the Lemmas \ref{linear-lem} and \ref{nonlin1} yield
$u(t)\rightharpoonup u_0$ in the weak$-\ast$ topology of $\dot{B}_{\infty,\infty}^{2/(\rho-1)}$ as $t\rightarrow 0^{+}$. The dependence of
the initial data can be obtained from Lemma \ref{linear-lem} and Lemma \ref{Aux-Lem}. Indeed, let $\bar{y}=L_{\alpha}(t)\bar{u_0}$ where $\bar{u_0}\in%
\mathcal{N}_{p,\mu,\infty}^{\sigma}$, then
\begin{align*}
\Vert u(t)-\bar{u}(t)\Vert _{X_q^p}&\leq \frac{1}{1-2^{\rho}K\varepsilon^{%
\rho-1}}\Vert L_{\alpha}(t)(u_0-\bar{u_0})\Vert_{X_q}\leq \frac{1}{%
1-2^{\rho}K\varepsilon^{\rho-1}}\Vert u_0-\bar{u_0}\Vert_{\mathcal{N}%
_{p,\mu,\infty}^{\sigma}}.
\end{align*}

\fin

\subsection{Proof of Theorem \protect\ref{teo2}}

The proof follows from analogous argument found in \cite[Theorem 3.3]%
{Almeida-Ferreira}. For the reader convenience, we indicate the main steps
of proofs.

\noindent \textbf{Item (i):} Let $M\in \mathcal{G}$ and $u_{0}$ be
antisymmetric, then $\Phi (x,t)=L_{\alpha }(t)u_{0}$ and $B_{\alpha }(u)$ is
antisymmetric. Indeed, in view of the orthogonality of $M$ and $%
u_{0}(Mx)=-u_{0}(x)$ we have
\begin{equation}
-\widehat{u_{0}}(\xi )=[u_{0}(M\cdot )]^{\wedge }(\xi )=\widehat{u_{0}}%
(M^{-1}\xi ),  \label{aux-sym1}
\end{equation}%
it follows that
\begin{align*}
\lbrack \Phi (Mx,t)]^{\wedge }(\xi )& =E_{\alpha }(-t^{\alpha }|M^{-1}\xi
|^{2})\widehat{a}(M^{-1}\xi ) \\
& =-E_{\alpha }(-t^{\alpha }|\xi |^{2})\widehat{a}(\xi ) \\
& =-\widehat{\Phi (x,t)}(\xi ),
\end{align*}%
this shows us that $L_{\alpha }(t)u_{0}$ is antisymmetric for each fixed $%
t>0.$ Similarly, we can show that $B_{\alpha }(u)$ is antisymmetric whether $%
u$ is also. So, employing an induction argument, one can prove that each
element $u_{k}$ of the Picard sequence
\begin{align}
u_{1}(x,t)& =\Phi (x,t)  \label{seq1} \\
u_{k}(x,t)& =\Phi (x,t)+B_{\alpha }(u_{k-1})(x,t),\;k=2,3,\cdots
\label{seq2}
\end{align}%
is antisymmetric. It follows that $u(x,t)$ is antisymmetric, for
all $t>0$. The symmetric property is analogous.

\noindent \textbf{Item (ii):} Let $\Phi (x,t)=L_{\alpha }(t)u_{0}$. In view
of $u_{0}$ be homogeneous of degree $-\frac{2}{\rho -1}$ we set $%
u_{0}(\lambda x)=\lambda ^{-2/(\rho -1)}u_{0}(x)$. It follows that
\begin{align*}
\lbrack \Phi (\lambda \cdot ,t)]^{\wedge }(\xi )& =E_{\alpha }(-t^{\alpha
}|\xi /\lambda |^{2})\widehat{u_{0}}(\xi /\lambda ) \\
& =\lambda ^{-\frac{2}{\rho -1}}\lambda ^{-n}E_{\alpha }(-t^{\alpha }|\xi
/\lambda |^{2})\widehat{u_{0}}(\xi ) \\
& =\lambda ^{-\frac{2}{\rho -1}}E_{\alpha }(-t^{\alpha }|\xi |^{2})\widehat{%
u_{0}}(\xi ) \\
& =\widehat{\lambda ^{-\frac{2}{\rho -1}}\Phi (\cdot ,t)}(\xi ),
\end{align*}%
that is, $\Phi (\lambda x,t)=\lambda ^{-\frac{2}{\rho -1}}\Phi (x,t)$. Now
proceeding like \textbf{Item (i)} we obtain that
\begin{equation*}
u(x,t)\equiv u_{\lambda }(x,t),\text{ for every }\lambda >0,
\end{equation*}%
in other words, $u$ is forward self-similar solution. \fin

\subsection{Proof of Theorem \protect\ref{teo3}}

We only show that (\ref{A2}) implies (\ref{A1}). The converse is left to the
reader. We have that

\begin{align}
t^{\eta }\Vert u(\cdot ,t)-v(\cdot ,t)\Vert _{\mathcal{M}_{q,\mu }}& \leq
t^{\eta }\Vert {L}_{\alpha }(t)(u_{0}-v_{0})\Vert _{\mathcal{M}_{q,\mu
}}+t^{\eta }\Vert B_{\alpha }(u)-B_{\alpha }(v)\Vert _{\mathcal{M}_{q,\mu }}
\notag \\
& :=t^{\eta }\Vert {L}_{\alpha }(t)(u_{0}-v_{0})\Vert _{\mathcal{M}_{q,\mu
}}+J_{1}(t)  \label{sa1}
\end{align}%
and
\begin{align}
\Vert u(\cdot ,t)-v(\cdot ,t)\Vert _{\mathcal{N}_{p,\mu ,\infty }^{\sigma
}}& \leq \Vert {L}_{\alpha }(t)(u_{0}-v_{0})\Vert _{\mathcal{N}_{p,\mu
,\infty }^{\sigma }}+\Vert B_{\alpha }(u)-B_{\alpha }(v)\Vert _{\mathcal{N}%
_{p,\mu ,\infty }^{\sigma }}  \notag \\
& \leq \Vert {L}_{\alpha }(t)(u_{0}-v_{0})\Vert _{\mathcal{N}_{p,\mu ,\infty
}^{\sigma }}+J_{2}(t).  \label{sa2}
\end{align}%
Using the inequality (\ref{nonlin-est1}), $\Vert u\Vert
_{X_{q}}\leq
2\varepsilon $ and $\left\Vert v\right\Vert _{X_{q}}\leq 2\varepsilon ,$ $%
J_{1}(t)$ can be estimated as
\begin{align}
J_{1}(t)& \leq Ct^{\eta }\int_{0}^{t}(t-s)^{\gamma
_{2}}\int_{0}^{s}r_{\alpha -1}(s-\tau )\Vert u(\tau )-v(\tau )\Vert _{%
\mathcal{M}_{q,\mu }}\left( \Vert u(\tau )\Vert _{\mathcal{M}_{q,\mu
}}^{\rho -1}+\Vert v(\tau )\Vert _{\mathcal{M}_{q,\mu }}^{\rho -1}\right)
d\tau ds  \notag \\
& \leq t^{\eta }2(2\varepsilon )^{\rho -1}C\int_{0}^{t}(t-s)^{\gamma
_{2}}\int_{0}^{s}r_{\alpha -1}(s-\tau )\tau ^{-\eta \rho }\Sigma _{1}(\tau
)d\tau ds,  \label{J1}
\end{align}%
where $\Sigma _{1}(\tau )=t^{\eta }\Vert u(\tau )-v(\tau )\Vert _{\mathcal{M}%
_{q,\mu }}$ and $\gamma _{2}=\eta \rho -\alpha -\eta .$ For $J_{2},$ in view
of (\ref{nonlin_B1-est11}), we have that
\begin{equation}
J_{2}(t)\leq (2^{\rho }\varepsilon ^{\rho -1})C\int_{0}^{t}(t-s)^{\gamma
_{1}}\int_{0}^{s}r_{\alpha -1}(s-\tau )\tau ^{-\eta \rho }\Sigma _{2}(\tau
)d\tau ds,  \label{J2}
\end{equation}%
where $\Sigma _{2}(\tau )=\Vert u(\tau )-v(\tau )\Vert _{\mathcal{N}_{p,\mu
,\infty }^{\sigma }}$ and $\gamma _{1}=\eta \rho -\alpha .$

Let us define $\Sigma (\tau )=\Sigma _{1}(\tau )+\Sigma _{2}(\tau ).$ After
performing a change of variables in (\ref{J1}) and (\ref{J2}), we get%
\begin{align}
J_{1}(t)+J_{2}(t)& \leq (2^{\rho }\varepsilon ^{\rho
-1})C\int_{0}^{1}(1-\sigma )^{\gamma _{2}}\sigma ^{\alpha -1-\eta \rho
}\int_{0}^{1}r_{\alpha -1}(1-z)z^{-\eta \rho }\Sigma (t\sigma z)dzd\sigma +
\notag \\
+& (2^{\rho }\varepsilon ^{\rho -1})C\int_{0}^{1}(1-\sigma )^{\gamma
_{1}}\sigma ^{\alpha -1-\eta \rho }\int_{0}^{1}r_{\alpha -1}(1-z)z^{-\eta
\rho }\Sigma (t\sigma z)dzd\sigma .  \label{aux-asymp1}
\end{align}
We claim that
\begin{equation}
\Pi :=\limsup_{t\rightarrow +\infty }\Sigma (t)=0,  \label{aux-asymp2}
\end{equation}%
which is enough for our purposes. For that, we take $\displaystyle%
\limsup_{t\rightarrow +\infty }$ in (\ref{aux-asymp1}) in order to obtain
{\footnotesize
\begin{align}
\limsup_{t\rightarrow +\infty }[J_{1}(t)+J_{2}(t)]& \leq (2^{\rho
}\varepsilon ^{\rho -1})C\int_{0}^{1}(1-\sigma )^{\gamma _{2}}\sigma
^{\alpha -1-\eta \rho }d\sigma \limsup_{t\rightarrow +\infty
}\int_{0}^{1}r_{\alpha -1}(1-z)z^{-\eta \rho }\Sigma (t\sigma z)dz  \notag \\
& +(2^{\rho }\varepsilon ^{\rho -1})C\int_{0}^{1}(1-\sigma )^{\gamma
_{1}}\sigma ^{\alpha -1-\eta \rho }d\sigma \limsup_{t\rightarrow +\infty
}\int_{0}^{1}r_{\alpha -1}(1-z)z^{-\eta \rho }\Sigma (t\sigma z)dz  \notag \\
& \leq (2^{\rho }\varepsilon ^{\rho -1})C\left( \int_{0}^{1}(1-\sigma
)^{\gamma _{2}}\sigma ^{\alpha -1-\beta \rho }d\sigma \int_{0}^{1}r_{\alpha
-1}(1-z)z^{-\eta \rho }dz\right) \Pi   \notag \\
& +(2^{\rho }\varepsilon ^{\rho -1})C\left( \int_{0}^{1}(1-\sigma )^{\gamma
_{1}}\sigma ^{\alpha -1-\eta \rho }d\sigma \int_{0}^{1}r_{\alpha
-1}(1-z)z^{-\eta \rho }dz\right) \Pi   \notag \\
& =(K_{1}+K_{2})(2^{\rho }\varepsilon ^{\rho -1})\Pi .  \label{aux-asymp3}
\end{align}
} It follows from (\ref{sa1}), (\ref{sa2}), (\ref{aux-asymp3}) and (\ref{A2}) that
\begin{align}
\Pi & \leq \limsup_{t\rightarrow +\infty }(t^{\eta }\Vert {L}_{\alpha
}(t)(u_{0}-v_{0})\Vert _{\mathcal{M}_{q,\mu }}+\Vert {L}_{\alpha
}(t)(u_{0}-v_{0})\Vert _{\mathcal{N}_{p,\mu ,\infty }^{\sigma
}})+\limsup_{t\rightarrow +\infty }[J_{1}(t)+J_{2}(t)]  \notag \\
& \leq 0+(K_{1}+K_{2})(2^{\rho }\varepsilon ^{\rho -1})\Pi =(2^{\rho
}\varepsilon ^{\rho -1}K)\Pi .  \label{aux-asymp4}
\end{align}%
So, due to $2^{\rho }\varepsilon ^{\rho -1}K<1,$ we get $\Pi =0$, as
required. \fin

\section*{Acknowledgments}

 M. de Almeida thanks to the Department of Mathematics of UNESP
- São José do Rio Preto for its kind hospitality during his visit.
M. de Almeida was partially supported by CNPQ, Grant
482428/2012-0, Brazil.

\end{document}